\documentclass[11pt]{amsart}
\usepackage{graphicx}
\usepackage{amssymb}
\usepackage{verbatim}
\usepackage{fullpage}

\usepackage{parskip}

\usepackage{amsthm}

\newtheorem{lemma}{Lemma}

\newtheorem{theorem}{Theorem}

\newtheorem{corollary}{Corollary}
\newtheorem{definition}{Definition}

\title{A characterization of triangulations of closed surfaces}
\author{J. Arocha, J. Bracho, N. Garc\'{i}a-Col\'{i}n, I. Hubard}

\begin{document}

%
%
%

\begin{abstract}
In this paper we prove that a finite triangulation of a connected closed surface is completely determined by its intersection matrix.
The \emph{intersection matrix} of a finite triangulation, $K$, is defined as $M_{K}=(dim(s_{i}\cap s_{j}))_{0\leq i,0\leq j}^{n-1}$, where $K_{2}=\{s_{0}, \ldots s_{n-1}\}$ is a labelling of the triangles of $K$.
\end{abstract}
\maketitle

\section{Introduction}

Within the theory of convex polytopes, the study of the combinatorial equivalence of $k$-skeleta of pairs polytopes which are not equivalent themselves has been of interest, this phenomena is referred to in the literature as ambiguity \cite{grunbaum1967convex}.

It is well known that for $k \geq \lfloor \frac{d}{2} \rfloor$ the $k$-skeleton of a convex polytope is not dimensionally ambiguous, this is, it defines the entire structure of its underlying $d$-polytope. However for $k < \lfloor \frac{d}{2} \rfloor$ the question is much more intricate.

One of the most interesting results in this direction is the solution to Perle's conjecture by P.Blind and R.Mani \cite{Blind1987} and, separately, by G. Kalai \cite{kalai1988simple} which states that the 1-skeleta of convex simple $d$-polytopes define their entire combinatorial structure. Or, on its dual version, that the dual graph (facet adjacency graph) of a convex simplicial d-polytope determines its entire combinatorial structure.

\section{Motivation \& contribution}

Allured by Perles' conjecture, we decided to explore the extent to which an adequate combination of combinatorial and topological assumptions would prove as powerful for characterising certain simplicial complexes. The purpose of this work is to present our first result, product of this exploration.

For topological assumption we will, in this instance, ask for the simplicial complex of study to be a connected closed surface. As for combinatorial assumption, one might be tempted to choose only to have the information provided by its dual graph. However, the dual graph of a triangulation of a closed surface does not provide enough information to characterise it, as there are some dual graphs to triangulations which have been shown  in \cite{mohar2004polyhedral} to have combinatorically different polyhedral embeddings.

Therefore, we will need to strengthen the combinatorial hypothesis. In order to do so we will introduce the concept of an intersection preserving mapping of simplices of a simplicial complex. 

\begin{definition} A bijective mapping $f: K_{d} \rightarrow K'_{d}$ between the sets of $d$-simplices of two simplicial complexes, $K$ and $K'$, is an intersection preserving mapping if for every pair of simplices $s, t \in K_{d}$ $dim(s \cap t)=dim(f(s)\cap f(t))$.
\end{definition}

Throughout this paper we will use the notation $K_{l}$ to refer to the set of $l$-dimensional simplices of the complex $K$.

Additionally, we will define two particular triangulations of the projective plane, which are of interest for this work.

\begin{definition} We define a 10-triangle triangulation of the projective plane, $T\mathbb{P}_{10}$, as the triangulation whose triangles have the vertex sets $(s_{i})_{0}=\{a_{i \mod 5}, a_{{i+1} \mod 5}, x\}$, and 
$(r_{i})_{0}=\{a_{i \mod 5}, a_{i+1 \mod 5}, a_{i-2 \mod 5}\}$ for $0 \leq i \leq 4$.
\end{definition}

\begin{definition} We define a 12-triangle triangulation of the projective plane, $T\mathbb{P}_{12}$, as the triangulation whose triangles have the vertex sets 
$(s_{i})_{0}=\{a_{i \mod 6}, a_{i+1 \mod 6}, x\}$, for $0 \leq i \leq 5$ and 
$(r_{i})_{0}=\{a_{i \mod 6}, a_{i+1 \mod 6}, a_{i+4 \mod 6}\}$ for $0 \leq i \leq 4$ even, and
$(r_{i})_{0}=\{a_{i \mod 6}, a_{i+1 \mod 6}, a_{i+3 \mod 6}\}$ for $0 \leq i \leq 4$ odd.
\end{definition}

\begin{figure}[htbp]
\begin{center}
\includegraphics[scale=0.5]{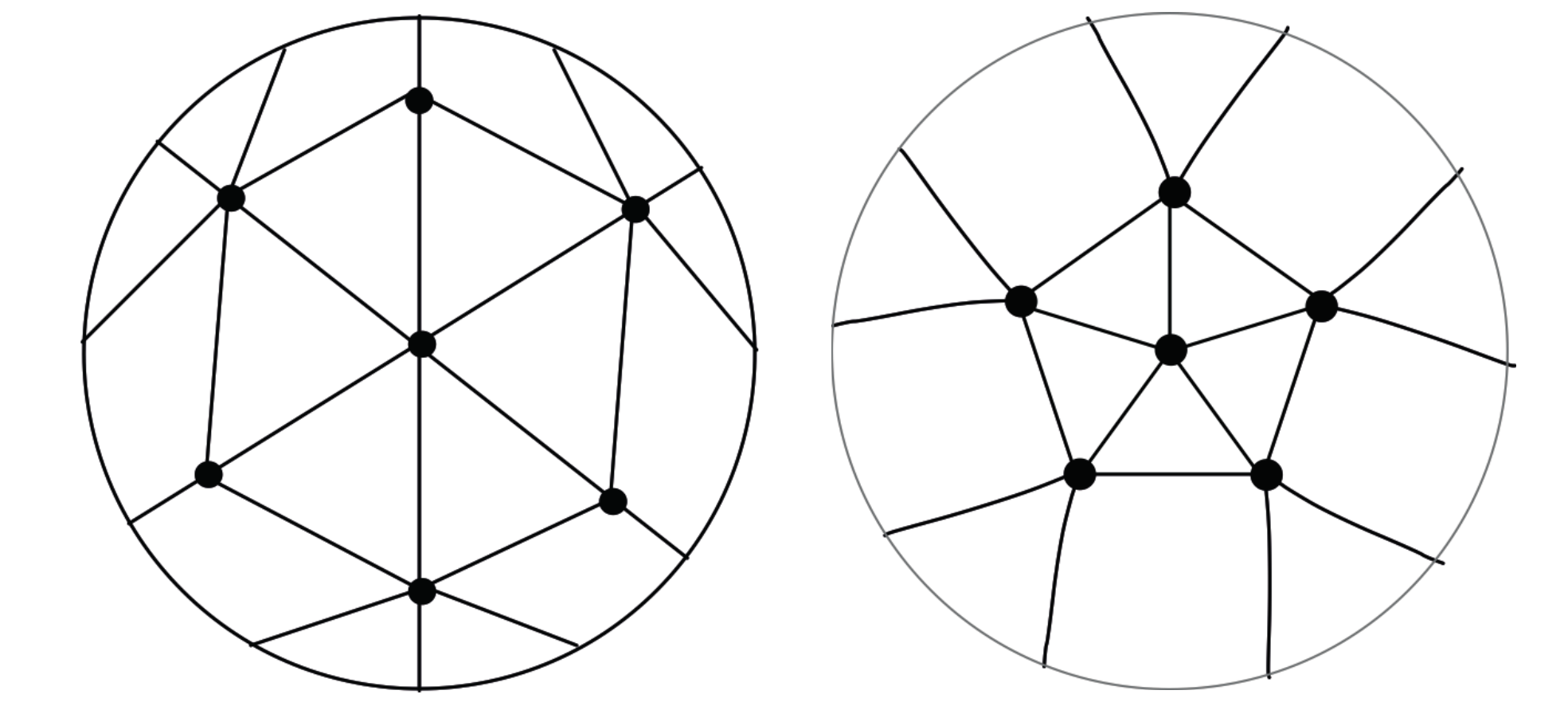}
\caption{Drawings for $T\mathbb{P}_{12}$ and $T\mathbb{P}_{10}$, respectively.}
\label{fig:only}
\end{center}
\end{figure}

Depictions of $T\mathbb{P}_{10}$ and $T\mathbb{P}_{12}$ are shown in figure \ref{fig:only}.

We now use the aforementioned definitions to state the main result:

\begin{theorem} Let $||K||$ and $||K'||$ be geometric realizations of finite triangulations which are homeomorphic to connected closed surfaces, and let $f: K_{2} \rightarrow K'_{2}$ be an intersection preserving mapping, then one of the following three statements holds:

\begin{enumerate}
\item $f$ can be extended into a bijective simplicial mapping between $K$ and $K'$
\item $f$ cannot be extended into a simplicial mapping between $K$ and $K'$, but both $||K||$ and $||K'||$ are $T\mathbb{P}_{10}$
\item $f$ cannot be extended into a simplicial mapping between $K$ and $K '$, but both $||K||$ and $||K'||$ are $T\mathbb{P}_{12}$.
\end{enumerate}
\end{theorem}

Consider the \emph{intersection matrix}, $M_{K}=(dim(s_{i}\cap s_{j}))_{0\leq i,0\leq j}^{n-1}$, of a  finite triangulation, $K$, where $K_{2}=\{s_{0}, \ldots s_{n-1}\}$ is a labelling of the triangles of $K$ then, in the spirit of Perles' conjecture, we can state the previous theorem as;

\begin{corollary}  A finite triangulation of a connected closed surface is completely determined by its intersection matrix. 
\end{corollary}

\section{Preliminaries}

One of the peculiarities of triangulations of a closed surface is that the neighbourhood of every vertex is a disk. Furthermore, the triangles incident to any vertex of such surface form the simplest of triangulations of a disk, namely an $n$-gon whose vertices are all linked by an edge to a central vertex in the centre of the $n$-gon. We start off by analysing the intersection patterns of such a structure.
 
\begin{definition} An \textbf{n-cycle} is the abstract triangulation $_{n}\Pi^{2}$ such that  $_{n}\Pi^{2}_{2}=\{s_0, s_1,\ldots, s_{n-1}\}$, $\dim(s_i \cap s_{i+1})=1$, $\dim(s_i \cap s_{j})=0$ for $|i-j|\geq2$ with $i,j \in \{0,\ldots, n-1\}$, and $\dim(s_{n-1} \cap s_{0})=1$.
\end{definition}

We will now focus on studying what other structures can have an intersection pattern equal to that of a triangulated disk. 

\begin{lemma} \label{lem:disk} The vertex sets of the triangles in an $n$-cycle, $_{n} \Pi^{2}_{2}=\{s_{0}, \ldots s_{n-1}\}$ can only take one of the following three types
\begin{enumerate}
\item $(s_{i})_{0}=\{a_{i \mod n}, a_{{i+1} \mod n}, x\}$ for all $0\leq i \leq n$, for any $n$;

\item $(s_{0})_{0}= \{a_{0}, a_{2}, a_{1}\}$, 
$(s_{1})_{0}= \{a_{1}, a_{3}, a_{2}\}$,
$(s_{2})_{0}= \{a_{2}, a_{4}, a_{3}\}$,
$(s_{3})_{0}= \{a_{3}, a_{0}, a_{4}\}$, and
$(s_{4})_{0}= \{a_{4}, a_{1}, a_{0}\}$, when $n=5$; or

\item $(s_{0})_{0}=\{a_{0}, a_{1}, a_{2}\}$, $(s_{1})_{0}= \{a_{1}, a_{2}, a_{4}\}$, 
$(s_{2})_{0}= \{a_{2}, a_{3}, a_{4}\}$,
$(s_{3})_{0}= \{a_{3}, a_{0}, a_{4}\}$,
$(s_{4})_{0}= \{a_{0}, a_{5}, a_{4}\}$, and
$(s_{5})_{0}= \{a_{5}, a_{2}, a_{0}\}$, when $n=6$.
\end{enumerate}
\end{lemma}

The proof of the lemma above consists of several parts and follows largely by a detailed analysis of the combinatorial structure of $n$-cycles of triangles.

It is easy to see that geometric realisations of the three types of triangulations associated to puzzles of $n$-cycles are an $n$-triangulation of a disk, a $5$-triangulation of a M\"{o}bius band and a $6$-triangulation of a M\"{o}bius band, respectively.

\section{Proof of the theorem}

We now proceed to proving the main theorem of this paper, using exhaustively the local and global implications of Lemma \ref{lem:disk}.
\begin{proof}

For each vertex $x \in K_{0}$ let $_{n_{x}}\Pi$ be the $n_{x}$-cycle around $x$, by hypothesis $||_{n_{x}}\Pi||$ is necessarily a disk. 

(1) If, for all $x \in K_{0}$, $||f(_{n_{x}}\Pi)||$ is also a disk, then the mapping $h: K_{0} \rightarrow K'_{0}$ such that $h(x)= \bigcap_{i=1}^{n_{x}} f(s_{i})$ is a bijective simplicial mapping.

Assume then that there is a vertex $x \in K_{0}$ such that $||f(_{n_{x}}\Pi)||$ is not a disk. 

(2) Suppose $||f(_{n_{x}}\Pi)||$ is the $5$-triangulation of the M\"{o}bius band described in Lemma \ref{lem:disk}.

Let $_{n_{x}}\Pi_{2}= \{ s_{0}, s_{1}, s_{2}, s_{3}, s_{4}\}$, where $(s_{i})_{0}=\{a_{i \mod 5}, a_{i+1 \mod 5}, x\}$ and  $f(_{n_{x}}\Pi)_{2}=\{ s_{0}', s_{1}', s_{2}', s_{3}', s_{4}'\}$ where 
$(s_{0}')_{0}= \{a_{0}', a_{2}', a_{1}'\}$,
$(s_{1}')_{0}= \{a_{1}', a_{3}', a_{2}'\}$,
$(s_{2}')_{0}= \{a_{2}', a_{4}', a_{3}'\}$,
$(s_{3}')_{0}= \{a_{3}', a_{0}', a_{4}'\}$, and
$(s_{4}')_{0}= \{a_{4}', a_{1}', a_{0}'\}$.

Given that $K'$ is also a closed surface, then each of the simplices $ s_{0}', s_{1}', s_{2}', s_{3}', s_{4}'$ has got a triangle adjacent to its remaining free edge. Let $r'_{i}$ be the simplices such that $dim(r'_{i}\cap s'_{i})=1$, then 
$(r_{0}')_{0}= \{a_{0}', a_{2}', x_{0}'\}$,
$(r_{1}')_{0}= \{a_{1}', a_{3}', x_{1}'\}$,
$(r_{2}')_{0}= \{a_{2}', a_{4}', x_{2}'\}$,
$(r_{3}')_{0}= \{a_{3}', a_{0}', x_{3}'\}$, and
$(r_{4}')_{0}= \{a_{4}', a_{1}', x_{4}'\}$. This  is
$(r_{i}')_{0}= \{a_{i\mod 5}', a_{i+2 \mod 5}', x_{i}'\}$. It follows that, $dim(r_{i}' \cap s_{j}')\geq 0$ for all $i\neq j$. 

Note that the interior of each of the edges $\{a_{i\mod 5}', a_{i+1 \mod 5}'\}$  is in the interior of the M\"{o}bius band, thus this edges cannot be repeated in any further simplex in the complex.

This implies that $x_{i} \not \in \{a_{0}', a_{1}', a_{2}', a_{3}', a_{4}'\}$, because, if this was the case, at least one of the edges $\{a_{i\mod 5}', a_{i+1 \mod 5}'\}$ would belong to $(r'_{i})_{1}$. Then, $dim(r_{i}' \cap s_{j}')= 0$

Let $r_{i}=f^{-1}(r'_{i})$, then $dim(r_{i} \cap s_{i})=1$ and $dim(r_{i} \cap s_{j})= 0$ for all $i\neq j$. As $(s_{i})_{0}=\{a_{i \mod 5}, a_{i+1 \mod 5}, x\}$ then $(r_{i})_{0}=\{a_{i \mod 5}, a_{i+1 \mod 5}, x_{i \mod 5}\}$.

Here $dim(r_{i\mod 5} \cap s_{i+1\mod 5})\geq 0$ and $dim(r_{i\mod 5} \cap s_{i-1\mod 5})\geq 0$ trivially, hence $a_{i-1 \mod 5}, a_{i+2 \mod 5} \not \in (r_{i})_{0}$. However, for $dim(r_{i\mod 5} \cap s_{i+2\mod 5})=0$ and $dim(r_{i\mod 5} \cap s_{i-2\mod 5})=0$ to be accomplished, necessarily $x_{i}=a_{i+3\mod 5}=a_{i-2\mod 5}$.

That is, $(r_{i})_{0}=\{a_{i \mod 5}, a_{i+1 \mod 5}, a_{i-2 \mod 5}\}$, hence the simplicial complex asociated to $\bigcup_{i=0}^{4}r_{i}$ is a $5$-triangulation of a M\'{o}bius band, where $dim(r_{i} \cap r_{i+2 \mod 5})=dim(r_{i} \cap r_{{i-2} \mod 5})=1$, and as necessarily $K_{2} = \bigcup_{i=0}^{4}r_{i} \cup \bigcup_{i=0}^{4}s_{i}$ and the geometric simplicial complexes associated to $\bigcup_{i=0}^{4}r_{i}$ and $\bigcup_{i=0}^{4}s_{i}$ are a M\"{o}bius band and a disk, respectively, then $||K||$ is equal to $T\mathbb{P}_{10}$.

The above also implies that $dim(r'_{i} \cap r'_{i+2 \mod 5})=dim(r'_{i} \cap r'_{i-2 \mod 5})=1$ then $v'=v'_{i}$ for all $i=1,\ldots 4$, so that $K'_{2} =\bigcup_{i=0}^{4}r'_{i} \cup \bigcup_{i=0}^{4}s'_{i}$, hence $||K'||$ is also equal to $T\mathbb{P}_{10}$.

(3) Suppose $||f(_{n_{x}}\Pi)||$ is the $6$-triangulation of the M\"{o}bius band described in Lemma \ref{lem:disk}.

Let $_{n_{x}}\Pi_{2}= \{ s_{0}, s_{1}, s_{2}, s_{3}, s_{4}, s_{5}\}$, where $(s_{i})_{0}=\{a_{i \mod 6}, a_{i+1 \mod 6}, x\}$ and  $f(_{n_{x}}\Pi)_{2}=\{ s_{0}', s_{1}', s_{2}', s_{3}', s_{4}', s_{5}'\}$ where 
$(s'_{0})_{0}=\{a'_{0}, a'_{1}, a'_{2}\}$,
$(s'_{1})_{0}= \{a'_{1}, a'_{2}, a'_{4}\}$,
$(s'_{2})_{0}= \{a'_{2}, a'_{3}, a'_{4}\}$,
$(s'_{3})_{0}= \{a'_{3}, a'_{0}, a'_{4}\}$,
$(s'_{4})_{0}= \{a'_{0}, a'_{5}, a'_{4}\}$, and
$(s'_{5})_{0}= \{a'_{5}, a'_{2}, a'_{0}\}$.

As $||K'||$ is a closed surface, then each of the simplices $\{ s_{0}', s_{1}', s_{2}', s_{3}', s_{4}', s_{5}'\}$ has got a triangle adjacent to its remaining free edge. Let $r'_{i}$ be the simplices such that $dim(r'_{i}\cap s'_{i})=1$, then 
$(r_{0}')_{0}= \{a_{0}', a_{1}', x_{0}'\}$,
$(r_{1}')_{0}= \{a_{1}', a_{4}', x_{1}'\}$,
$(r_{2}')_{0}= \{a_{2}', a_{3}', x_{2}'\}$,
$(r_{3}')_{0}= \{a_{3}', a_{0}', x_{3}'\}$,
$(r_{4}')_{0}= \{a_{4}', a_{5}', x_{4}'\}$, and
$(r_{5}')_{0}= \{a_{2}', a_{5}', x_{5}'\}$.

Here it follows that, $dim(r_{i}' \cap s_{j}')\geq 0$ for all $i\neq j$, except for the pairs $i=0 \text{ and } j=2$, $i=1\text{ and }j=5$, $i=2\text{ and }j=4$, $i=3\text{ and }j=1$, $i=4\text{ and }j=0$ and $i=5\text{ and }j=3$; for these exceptions the intersection might be empty.

The above implies that, if $r_{i}=f^{-1}(r'_{i})$, then $dim(r_{i} \cap s_{i})=1$ and $dim(r_{i} \cap s_{j})\geq 0$ for all $i\neq j$, except for the pairs $i=0 \text{ and } j=2$, $i=1\text{ and }j=5$, $i=2\text{ and }j=4$, $i=3\text{ and }j=1$, $i=4\text{ and }j=0$, and $i=5\text{ and }j=3$; for these exceptions the intersection might be empty.

As $(s_{i})_{0}=\{a_{i \mod 6}, a_{i+1 \mod 6}, x\}$ then the vertex sets of the $r_{i}$'s are $(r_{i})_{0}=\{a_{i \mod 6}, a_{i+1 \mod 6}, x_{i \mod 6}\}$.

Note that $x_{0}' \not \in \{a'_{0}, a'_{1}, a'_{2}, a'_{4}, a'_{5}\}$ as the edges $\{a'_{0}, a_{2}\}$, $\{a'_{0}, a_{4}\}$, $\{a'_{0}, a_{5}\}$, $\{a'_{0}, a_{4}\}$ are edges whose interior is in the interior of the M\"{o}bius band. Thus we might have $v'_{0}=a'_{3}$, however if that was the case
$dim(r_{0}' \cap r_{3}')= 1$, and $dim(r_{0} \cap r_{3})= 1$, but this is not possible. Then necessarily  $x_{0}' \not \in \{a'_{0}, a'_{1}, a'_{2}, a'_{3}, a'_{4}, a'_{5}\}$. 

Using an argument analogous to the one in the previous case, we deduce that for each $i$, $x_{i}' \not \in \{a'_{0}, a'_{1}, a'_{2}, a'_{3}, a'_{4}, a'_{5}\}$; so that $dim(r_{i}' \cap s_{j}')= 0$ for all $i\neq j$, except for the pairs $i=0\text{ and }j=2$, $i=1\text{ and }j=5$, $i=2\text{ and }j=4$, $i=3\text{ and }j=1$, $i=4\text{ and }j=0$, and $i=5\text{ and }j=3$, for which the intersection is empty. 

The above implies $dim(r_{i} \cap s_{i})=1$ and $dim(r_{i} \cap s_{j})= 0$ for all $i\neq j$, except for the pairs $i=0\text{ and } j=2$, $i=1\text{ and } j=5$, $i=2\text{ and } j=4$, $i=3\text{ and } j=1$, $i=4\text{ and } j=0$, and $i=5\text{ and } j=3$, for which the intersection is empty. Hence, in order to accomplish the intersection dimensions indicated by the puzzle necessarily, $x_{0}=x_{1}= a_{4}$, $x_{2}=x_{3}= a_{0}$, $x_{4}=x_{5}= a_{2}$, thus;
$r_{0}=\{a_{0}, a_{1}, a_{4}\}$,
$r_{1}=\{a_{1}, a_{2}, a_{4}\}$,
$r_{2}=\{a_{2}, a_{3}, a_{0}\}$,
$r_{3}=\{a_{3}, a_{4}, a_{0}\}$,
$r_{4}=\{a_{4}, a_{5}, a_{2}\}$, and
$r_{5}=\{a_{0}, a_{5}, a_{2}\}$.

Therefore, the simplicial complex associated to $\bigcup_{i=0}^{5} r_{i}$ is a $6$-triangulation of a M\"{o}bius band and $K_{2}= \bigcup_{i=0}^{5} s_{i} \cup \bigcup_{i=0}^{5} r_{i}$, so that $||K||$ is equal to $T\mathbb{P}_{12}$.

The implication for $K'$ is that $dim(r'_{0} \cap r'_{1})=1,\; dim(r'_{1} \cap r'_{4})=1, \; dim(r'_{4} \cap r'_{5})=1, \; dim(r'_{5} \cap r'_{2})=1, \; dim(r'_{2} \cap r'_{3})=1, \; dim(r'_{3} \cap r'_{0})=1$, which in turn implies $v'=v'_{i}$ for all $i \in \{0,\ldots 5\}$ and, further, $K'_{2}= \bigcup_{i=0}^{5} s'_{i} \cup \bigcup_{i=0}^{5} r'_{i}$, so that $K'$ is also equal to $T\mathbb{P}_{12}$.
\end{proof}

\bibliographystyle{plain}
\bibliography{simplicialpuzzles}

\end{document}